%% file: main.tex
\documentclass[letterpaper, 10 pt, conference]{ieeeconf}
\IEEEoverridecommandlockouts
\usepackage{cite}
\usepackage{mathrsfs}
\usepackage{algorithmic}
\usepackage{graphicx}
\usepackage{textcomp}
\usepackage{xcolor}
\usepackage{svg}
\usepackage{float}
\usepackage{subcaption}
\usepackage{caption}
\usepackage{titlesec}
\usepackage[matmatrix, opt, reviewmark]{icpslab}
\usepackage{etoolbox}
\usepackage{booktabs}
\usepackage{multirow}
\usepackage{siunitx}
\usepackage{hyperref}
\usepackage[table]{xcolor}
\let\labelindent\relax
\usepackage{enumitem}

\newtheorem{assumption}{Assumption}

\newcommand{\SBMPC}{SBMPC\xspace}

\begin{document}

\title{\LARGE \bf
Learning-enabled Acceleration of Scenario-based \\ Model Predictive Control
}
\author{Trinh Tran, Binh Nguyen, Truong X. Nghiem
\thanks{*This material is based upon work supported by the U.S.\ National Science Foundation under Grants No.~2449927 and No.~2514584.}%
\thanks{All authors are with Department of Electrical and Computer Engineering,
University of Central Florida, Orlando, FL 32816, USA. \tt\small \{tr586225, thanhbinh.nguyen, truong.nghiem\}@ucf.edu
}
}

\maketitle

\begin{abstract}
Scenario-based model predictive control (\SBMPC) is a variant of model predictive control (MPC) that explicitly accounts for uncertainty by optimizing control actions over multiple predicted scenarios.
However, its computational complexity increases rapidly with the number of scenarios and prediction horizon, limiting its applicability to real-time planning and control.
This paper presents a learning-accelerated Alternating Direction Method of Multipliers (ADMM) algorithm for efficiently solving \SBMPC  problems
by leveraging parallel computing and Moreau envelope learning, while maintaining high solution accuracy.
We reformulate the \SBMPC  problems into consensus forms that can be  decomposed via ADMM, separating the scenario-dependent dynamics from non-anticipativity constraints and enabling parallel updates across scenarios and time steps. 
Building on this decomposition, we utilize a learning-to-optimize scheme that leverages Moreau envelope learning of the cost function to accelerate the primal update in ADMM, thereby reducing computation time.
The proposed framework is evaluated on a microgrid energy management problem subject to load and renewable generation uncertainties.
Comparisons with IPOPT and MadNLP, two popular and modern nonlinear programming solvers, demonstrate substantial computational speedups while maintaining reliable closed-loop control performance.
\end{abstract}

\input{intro}
\section{Problem Formulation} \label{sec:problem-formulation}
We consider a discrete-time linear system
\begin{equation}
    x_{t+1} = A x_t + B u_t + d_t,
\end{equation}
where $x_t \in \mathbb{R}^n$ denotes the system state vector, 
$u_t \in \mathbb{R}^m$ is the control input vector,
$d_t \in \mathbb{R}^n$ represents an unmeasured disturbance vector with a possibly unknown distribution $\mathbb{D}_d$,
and $A\in \mathbb{R}^{n \times n}, B \in \mathbb{R}^{n \times m}$ are the system matrices.

To account for disturbance in a tractable way, the possible future disturbance is approximated by a finite set of scenarios
$s = 1,...,S$. 
Each scenario $s$ captures a realization of the uncertain model parameters and disturbance sequence over the prediction horizon.
We use $x_{s,k}$ and $u_{s,k}$ to denote the state and input predicted $k$ steps ahead under scenario $s$ at time $t$.
Consider the following scenario MPC problem, where the disturbance is sampled from its distribution $\mathbb{D}_d$ to generate $S$ scenarios of the disturbance sequence $\{d_{s,t+k}\}_{k=0}^{N-1}$:

\begin{subequations}
\label{eq:scenario-MPC}
\begin{align}
\minimize\;\; & 
    \frac{1}{S} \sum_{s=1}^{S} \left(\sum_{k=0}^{N-1} c_k(x_{s,k}, u_{s,k}) + c_N(x_{s,N}) \right)\\
    \text{s.t.}\;\;
    & x_{s,k + 1} = A x_{s,k} + B u_{s,k} + d_{s,t+k},\label{eq:scenario-MPC:dynamics}\\
    &x_{s,0} = x_t, \;\forall s=1, \cdots, S,\\
    &(x_{s,k}, u_{s,k}) \in \mathcal{X} \times \mathcal{U}, \label{eq:scenario-MPC:state-control}\\
    &\qquad \quad \forall k = 0, \cdots, N-1, \, s=1, \cdots, S, \nonumber\\
    &P \mathbf{u} = 0\text, \label{eq:scenario-MPC:non-anticipativity}
\end{align}
\end{subequations}
where $x_t$ is the current state, 
$N$ is the prediction horizon,
$\mathbf{x}_s = (x_{s,1},\dots,x_{s,N}) \in \mathbb{R}^{Nn}$,
$\mathbf{u}_s = (u_{s,0},\dots,u_{s,N-1}) \in \mathbb{R}^{Nm}$, 
$c_k(\cdot, \cdot)$ and $c_N(\cdot)$ are convex, possibly nonsmooth, stage and terminal cost functions,
and $\mathcal{X}$ and $\mathcal{U}$ are the state and input constraint sets.
If the matrix $P$ is defined as
\begin{equation*}
  P = \mm{\IdentityMatrix_{N m}, -\IdentityMatrix_{N m}, 0, 0, \dots, 0;
        0, \IdentityMatrix_{N m}, -\IdentityMatrix_{N m}, 0, \dots, 0;
        \dots;
        0, 0, \dots, 0, \IdentityMatrix_{N m}, -\IdentityMatrix_{N m}} \text,
\end{equation*}
then the non-anticipativity constraint~\eqref{eq:scenario-MPC:non-anticipativity} enforces that $u_{1,k} = u_{2,k} = \dots = u_{S,k}$ for $k = 0,\dots,N-1$, \ie the control input is scenario-independent.
On the other hand, if 
\begin{equation*}
  P = \mm{P_0, -P_0, 0, 0, \dots, 0;
        0, P_0, -P_0, 0, \dots, 0;
        \dots;
        0, 0, \dots, 0, P_0, -P_0}\text,
\end{equation*}
where $P_0 = \mm{\IdentityMatrix_{n_u}, 0, \dots, 0}$, then the constraint~\eqref{eq:scenario-MPC:non-anticipativity} enforces that $u_{1,0} = u_{2,0} = \dots = u_{S,0}$, \ie the control input at the first step is scenario-independent.

The following standard assumptions are made.
\begin{assumption}\label{ass:convexity}
(i) The full system state is available 
at each time step $t$.
(ii) The control input constraint set $\mathcal{U}$ is bounded and convex.
(iii) The state constraint set $\mathcal{X}$ is convex.
(iv) Convexity holds for the stage $c_k$ and terminal $c_N$ costs.
\end{assumption}

The dynamics constraint for each scenario $s$ can be rewritten as
$\mathbf{x}_s = M_A \mathbf{x}_s + M_B \mathbf{u}_s + M_0 x_t + \mathbf{d}_s$, or equivalently
\begin{equation}
\label{eq:rewrite-dynamics}
    M_G\mm{\mathbf{x}_s; \mathbf{u}_s} = -M_0 x_t - \mathbf{d}_s,
\end{equation}
where $\mathbf{d}_s = (d_{s,t},\dots,d_{s,t+N-1}) \in \mathbb{R}^{Nn}$, 
$M_G = \mm{M_A - \IdentityMatrix_{Nn}, M_B} \in \mathbb{R}^{Nn \times (Nn+Nm)}$, with
$M_A \in \mathbb{R}^{Nn \times Nn}$, $M_B \in \mathbb{R}^{Nn \times Nm}$, and $M_0 \in \mathbb{R}^{Nn \times n}$ defined as
\begin{align*}
M_A &= \begin{bmatrix} 
0 &  &  &  \\ 
A & 0 &  &  \\ 
 & \ddots & \ddots &  \\ 
 &  & A & 0 
\end{bmatrix}, \;
M_B = \IdentityMatrix_N \otimes B, \;
M_0 = \begin{bmatrix} 
A \\ 0 \\ \vdots \\ 0 
\end{bmatrix}.
\end{align*}
Therefore, under Assumption~\ref{ass:convexity}, the dynamics constraint set 
\[
\mathcal{C}_s := \left\{ (\mathbf{x}_s, \mathbf{u}_s) \in \mathcal{X}^N \times \mathcal{U}^N \;\middle|\; M_G \begin{bmatrix} \mathbf{x}_s\\ \mathbf{u}_s \end{bmatrix} = -M_0 x_t - \mathbf{d}_s \right\}
\]
is convex and parameterized by $x_t$ and $\mathbf d_s$.

Although the scenario-based formulation~\eqref{eq:scenario-MPC} provides a deterministic approximation for the stochastic MPC problem, solving it incurs significant computational overhead. 
This is because the size of the problem grows quickly with the number of scenarios, prediction horizon, and system dimensions.
Furthermore, the coupling constraints~\eqref{eq:scenario-MPC:non-anticipativity} and \eqref{eq:rewrite-dynamics} require solving a full trajectory optimization for each scenario. 
This motivates a decomposition approach that splits across both scenarios and time steps to enable parallel computation.

\section{Learning-acclerated ADMM for SBMPC} \label{sec:proposed-method}

\subsection{ADMM algorithm for scenario-based MPC}

We assume that all scenarios share the same stage and terminal cost functions, while the disturbance realizations differ across scenarios. 
We remove $\frac{1}{S}$ from the objective function for simplicity and rewrite~\eqref{eq:scenario-MPC} in consensus form as
\begin{subequations}
\label{eq:consensus-scenario-MPC}
\begin{align}
\minimize\;&\sum_{s=1}^{S}\!
\left(\sum_{k=0}^{N-1}\!c_k(x_{s,k},\!u_{s,k}) \!+\! c_N(x_{s,N}) \!+\! {I}_{\mathcal{C}_s} (\mathbf{y}_{s},\!\mathbf{v}_{s}) \right)
 \nonumber\\
& \qquad\qquad +{I}_{\mathcal{P}}(\mathbf{p}_{1},\dots,\mathbf{p}_{S})\label{eq:consensus-scenario-MPC:obj} \\
\text{s.t.}\;&\mathbf{x}_{s} = \mathbf{y}_{s},\; s = 1,\dots,S\text,\\
&\mathbf{u}_{s} = \mathbf{p}_{s},\; s = 1,\dots,S\text,\\
&\mathbf{v}_{s} = \mathbf{p}_{s},\; s = 1,\dots,S\text,
\end{align}
\end{subequations}
where $\mathbf{y}_s \in \mathbb{R}^{Nn}$, $\mathbf{v}_s \in \mathbb{R}^{Nm}$, and $\mathbf{p}_s \in \mathbb{R}^{Nm}$ are auxiliary variables,
and $\mathcal{P}=\{\mathbf{p}: P\mathbf{p}=0\}$, $\mathbf{p}=(\mathbf{p}_1,\dots,\mathbf{p}_S),\, \mathbf{p}_s\in\mathbb{R}^{Nm}$.
Here, ${I}_{\mathcal{C}_s} (\mathbf{y}_{s}, \mathbf{v}_{s})$ is the indicator function for the dynamics constraint set $\calC_s$, which equals $0$ if $(\mathbf{y}_{s}, \mathbf{v}_{s})\in \mathcal{C}_s$,and $+\infty$ otherwise.
Similarly, ${I}_{\mathcal{P}}(\cdot)$ is the indicator function for the set $\mathcal{P}$.

For the Lagrange multipliers $\lambda_{y,s}$ associated with $\mathbf{x}_{s} -\mathbf{y}_{s} = 0$,
$\lambda_{v,s}$ associated with $\mathbf{v}_{s} - \mathbf{p}_{s} = 0$,
and $\lambda_{p,s}$ associated with $\mathbf{u}_{s} - \mathbf{p}_{s} = 0$,
and for a penalty parameter $\rho > 0$,
the augmented Lagrangian is
\begin{align*}
  \mathcal{L}_{\rho}&(\mathbf{x}, \mathbf{u}, \mathbf{y}, \mathbf{v}, \mathbf{p}, \lambda_y, \lambda_v, \lambda_p) = 
  \sum_{s=1}^{S}\sum_{k=0}^{N-1}c_k(x_{s,k}, u_{s,k}) \\
  &+ \sum_{s=1}^{S}c_N(x_{s,N}) +\sum_{s=1}^{S} \mathcal{I}_{\mathcal{C}_s}(\mathbf{y}_{s}, \mathbf{v}_{s})
  + \mathcal{I}_{\mathcal{P}}(\mathbf{p}_{1}, \dots, \mathbf{p}_{S}) \\
  &+ \sum_{s=1}^{S} \left( \lambda_{y,s}^{\top}(\mathbf{x}_{s} - \mathbf{y}_{s})
  + \lambda_{v,s}^{\top}(\mathbf{v}_{s} - \mathbf{p}_{s})
  + \lambda_{p,s}^{\top}(\mathbf{u}_{s} - \mathbf{p}_{s}) \right) \\
  &+ \frac{\rho}{2}\sum_{s=1}^{S} \left( \|\mathbf{x}_{s} - \mathbf{y}_{s}\|^2_2
  + \|\mathbf{v}_{s} - \mathbf{p}_{s}\|^2_2
  + \|\mathbf{u}_{s} - \mathbf{p}_{s}\|^2_2 \right).
\end{align*}
The ADMM iterates for the scenario MPC problem~\eqref{eq:consensus-scenario-MPC} are
\begin{subequations}
\label{eq:admm-consensus-scenario-MPC}
\begin{align}
  (x_{s,k}^{i+1},u_{s,k}^{i+1})
  &=\operatorname*{arg\,min}_{x_{s,k},u_{s,k}}\Big\{c_k(x_{s,k},u_{s,k}) \nonumber\\
  &+\frac{\rho}{2}\|x_{s,k}-y_{s,k}^i+\rho^{-1}\lambda_{y,s,k}^i\|^2_2
  \nonumber\\
  &+\frac{\rho}{2}\|u_{s,k}-p_{s,k}^i+\rho^{-1}\lambda_{p,s,k}^i\|^2_2
  \Big\}, \nonumber\\
  &\qquad\;s = 1,\dots,S;\; k = 0,\dots,N-1\text,\label{eq:admm-consensus-scenario-MPC:primal}\\
  x_{s,N}^{i+1} &= \arg\min_{x_{s,N}} \Big\{ c_N(x_{s,N}) \nonumber \\
  &+ \frac{\rho}{2} \left\| x_{s,N} - y_{s,N}^i + \rho^{-1} \lambda_{y,s,N}^i \right\|_2^2 \Big\},\label{eq:admm-consensus-scenario-MPC:terminal_primal}
  \\(\mathbf{y}_{s}^{i+1},\mathbf{v}_{s}^{i+1}) &=\Pi_{\mathcal{C}_s}\left(\mm{\mathbf{x}_{s}^{i+1}+\lambda_{y,s}^i; \mathbf{p}_{s}^i \!-\! \lambda_{v,s}^i}\right),\;\; s = 1,\dots,S\text,
  \label{eq:admm-consensus-scenario-MPC:consensus}\\
  \mathbf{p}^{i+1}
  &=\Pi_{\mathcal{P}}\left(\frac{1}{2}\left(\mathbf{u}^{i+1} \!+\!\mathbf{v}^{i+1}\right)
  \!+\! \frac{1}{2\rho}\left(\lambda_{p}^i+\lambda_{v}^i \!\right)\!\!
  \right)\!\text,\label{eq:admm-consensus-scenario-MPC:non-anticipative}\\
  \lambda^{i+1} &= \lambda^i + \rho\left(\mm{\mathbf{x}^{i+1}; \mathbf{v}^{i+1}; \mathbf{u}^{i+1}}-\mm{\mathbf{y}^{i+1}; \mathbf{p}^{i+1}; \mathbf{p}^{i+1}}\right)\text,
  \label{eq:admm-consensus-scenario-MPC:lambda}
\end{align}
\end{subequations}
where $\lambda = [\lambda_y^\top, \lambda_v^\top, \lambda_p^\top]^\top$ is the lumped Lagrange multiplier vector,
$\Pi_{\mathcal{C}_s}(\cdot)$ and $\Pi_{\mathcal{P}}(\cdot)$ are the Euclidean projection operators onto sets $\mathcal{C}_s$ and ${\mathcal{P}}$, respectively.

These subproblems~\eqref{eq:admm-consensus-scenario-MPC:primal} are separable in $s$ and $k$, while the terminal updates~\eqref{eq:admm-consensus-scenario-MPC:terminal_primal} are separable across scenarios.
Therefore, they can be solved in parallel across scenarios and time steps.
However, their computational cost can still be significant when the number of scenarios or the prediction horizon becomes large. 

\subsection{Moreau Envelope Learning for ADMM}
In this part, we apply the Learning-Enabled ADMM Framework developed in our previous work~\cite{nguyen2026leaflearningenabledadmmframework} to accelerate the primal update step.
Rather than directly solving the optimization problems in ~\eqref{eq:admm-consensus-scenario-MPC:primal} and \eqref{eq:admm-consensus-scenario-MPC:terminal_primal}, LEAF leverages the Moreau envelope to provide a smooth, learnable approximation of the primal step.

The primal update can be equivalently written as the proximal operator of the local objective function as
\begin{equation*}
    z_{s,k}^{i+1} = \text{prox}_{\rho^{-1} c_k} (w_{s,k}^i),\quad k=0, \dots, N, \; s=1, \dots, S
\end{equation*}
where $z_{s,k}^{i+1}= [x_{s,k}^{i+1},u_{s,k}^{i+1}]^\top$, $w_{s,k}^i = [y_{s,k}^i, p_{s,k}^i]^\top - \rho^{-1} [\lambda_{y,s,k}^i, \lambda_{p,s,k}^i]^\top$.
Note that $u_{s,N} = p_{s,N} = \lambda_{p,s,N} = 0_m,\; s = 1, \dots, S.$
For more details on proximal operators and their role in convex optimization, we refer readers to~\cite{neal2014proximal}.

The corresponding Moreau envelope~\cite{bauschke2020correction} of $c_k$ with parameter $\rho^{-1}$ is defined as
\begin{equation*}
M_{\rho^{-1}c_k}(w_{s,k}^i)=\min\left\{c_k(z_{s,k})+ \frac{\rho}{2}\left\|z_{s,k} - w_{s,k}^i\right\|^2_2\right\}.
\end{equation*}
The proximal operator and the Moreau envelope satisfy the following relationship
\begin{equation}
\label{eq:prox-moreau-relationship}
\text{prox}_{\rho^{-1} f}(w)=w-\rho^{-1} \nabla M_{\rho^{-1} f}(w).
\end{equation}
Eventually, the primal update can be equivalently written as
\begin{equation*}
    \label{eq:primal-update-by-ME}
 z_{s,k}^{i+1} = w_{s,k}^i - \rho^{-1}\nabla M_{\rho^{-1}c_k}(w_{s,k}^i).
\end{equation*}
The Moreau envelope preserves the convexity of $c_k$ while providing a differentiable approximation, even when $c_k$ is nonsmooth.
Therefore, instead of directly learning the solutions of the nonsmooth optimization problem~\eqref{eq:admm-consensus-scenario-MPC:primal},\;\eqref{eq:admm-consensus-scenario-MPC:terminal_primal}, our learning-based method approximates the gradient of the Moreau envelope and uses it to recover the primal update.

For any observed input $w_{s,k}^i$, we train an input-convex neural network~\cite{pmlr-v70-amos17b} model $\mathcal{N}_{\theta,k}$ with weights $\theta_k$ to approximate the Moreau envelope $M_{\rho^{-1}c_k}(w_{s,k}^i)$.
Equivalently, the neural network outputs $\hat{M}_{\theta,k}(w_{s,k}^i)=\mathcal{N}_{\theta,k}(w_{s,k}^i)$, and from that we can obtain its gradient $\nabla \hat{M}_{\theta,k}(w_{s,k}^i)$ according to~\eqref{eq:prox-moreau-relationship}.
The gradient of the learned envelope is then used to reconstruct a primal update as follows
\begin{equation*}
\hat z_{s,k}^{i+1}=w_{s,k}^i-\rho^{-1}\nabla \hat M_{\theta,k}(w_{s,k}^i),
\label{eq:learned-primal-update}
\end{equation*}
where $\hat z_{s,k}^{i+1}=[\hat x_{s,k}^{i+1},\hat u_{s,k}^{i+1}]^\top$ with $\hat{x}_{s,k}^{i+1}$ and $\hat{u}_{s,k}^{i+1}$ are the approximated primal variables.
Thus, the expensive solution of the local unconstrained optimization problem in~\eqref{eq:admm-consensus-scenario-MPC:primal} is replaced by a forward evaluation of the neural network and a gradient computation.
Once trained, the learned update can be reused across ADMM iterations and across scenario-time blocks that share the same local cost structure.

The above training goal could be realized by minimizing the following loss function for the neural network $\mathcal{N}_{\theta,k}$ in a supervised-learning manner:
\begin{equation}
\label{eq:mel-training-loss}
\begin{aligned}
L_k(\theta_k) =&
\sum_{j=1}^{n_D}\Big(
\left\|\hat M_{\theta,k}(w_j)-M_{\rho^{-1}c_k}(w_j)\right\|_2^2\\
&+ \alpha\left\|\nabla \hat M_{\theta,k}(w_j)-\nabla M_{\rho^{-1}c_k}(w_j)\right\|_2^2\\
&+ \beta\left\|\max\left(0,\hat M_{\theta,k}(w_j)-c_k(w_j)\right)\right\|_2^2
\Big),
\end{aligned}
\end{equation}
where $\alpha>0$ and $\beta>0$ are weighting coefficients, $n_D$ is the number of supervised training samples, and $j$ indexes each training sample.
The last term penalizes violations of the Moreau-envelope property $M_{\lambda f}(x)\leq f(x)$ for all $x$, by enforcing $\hat M_{\theta,k}(w_j)\leq c_k(w_j)$ on the training samples.
The scenario index $s$ does not appear explicitly in~\eqref{eq:mel-training-loss} because the network is trained on local inputs $w_j$ drawn from the collection of inputs $w_{s,k}^i$ across scenarios, time steps, and ADMM iterations.
After collection, each $w_j$ is treated as one sample for the learning map associated with $c_k$.

\section{Experimental results}\label{sec:experimental-results}
In this section, we evaluate the proposed method on the optimal control task of energy management in microgrids.
\subsection{System description}
We consider the optimal energy management of a microgrid with renewable energy integration and battery energy storage system (BESS)~\cite{cortes-aguirreEconomicMPC2025}.
The system consists of controllable loads, uncontrollable loads, a BESS, a renewable generator, and a utility-grid connection.
The objective is to schedule the BESS charging/discharging power to maximize the average profit.

The uncertainty in the grid arises from fluctuations in renewable energy generation and load demand.
We represent this uncertainty by a finite scenario set $\mathcal{S} = \{1,\ldots,S\}$.
We formulate the operation problem for the microgrid over the planning horizon $N$ with time step $\Delta t$.
For each scenario $s \in \mathcal{S}$, the net uncontrollable power disturbance at step $k$
\begin{equation}
\label{eq:em-disturbance}
  d_{s,k} = P^u_{s,k} - P^R_{s,k}\text,
\end{equation}
where $P^u_{s,k} \ge 0$ is the uncontrollable load and $P^R_{s,k} \ge 0$ is the renewable generation.
A positive $d_{s,k}$ means uncontrollable demand exceeds renewable supply.
The state is the BESS state of charge (SOC) $q_{s,k}$,
and the control $u_k = [P^c_k, P^B_k]^\top$is the scenario-independent schedule.
$P^c_k$ is the controllable load and $P^B_k$ is the BESS charging/discharging power.
A positive $P^B_k$ represents discharging and a negative value represents charging.
The SOC dynamics for scenario $s$ are
\begin{equation}
\label{eq:em-SOC-dynamics}
  q_{s,k+1} = q_{s,k} - \frac{\Delta t}{b}\,P^B_k\text,
\end{equation}
The schedule is scenario-independent because decisions must be made before the uncertainty is realized (non-anticipativity).
Note that the disturbance $d_{s,k}$ does not appear in~\eqref{eq:em-SOC-dynamics}.
Instead, it enters through the power balance $P^c_k + P^u_{s,k} = D_{s,k} + P^B_k + P^R_{s,k}$,
which gives the net power drawn from the utility grid:
\begin{equation}
\label{eq:em-grid-exchange}
  D_{s,k} = P^c_k + d_{s,k} - P^B_k\text.
\end{equation}

The stage cost for scenario $s$ at step $k$ is the sum of an energy cost, a peak-demand cost, and a discomfort cost:
\begin{equation}
\label{eq:em-stage-cost}
\begin{aligned}
  c_k(q_{s,k},\, u_k,\, d_{s,k})
  &= p_k \Delta t \!\left(D_{s,k} + \frac{1-\mu}{2\sqrt{\mu}}\,|P^B_k|\right)\\
  &+p_p \max(D_{s,k},\, 0)
  +\eta\, g(P^c_k)\text,
\end{aligned}
\end{equation}
where $p_k$ is the electricity price,
$\mu \in (0,1]$ is the round-trip efficiency of the BESS,
$p_p \ge 0$ weights the peak-demand charge,
and $\eta \ge 0$ weights the discomfort penalty.
The term $\frac{1-\mu}{2\sqrt{\mu}}|P^B_k|$ accounts for the energy lost due to BESS cycling.
The convex discomfort function $g(\cdot)$ penalizes consumption below the preferred level $\gamma$:
\begin{equation*}
  g(P) =
  \begin{cases}
    +\infty, & P \le 0, \\
    \frac{\gamma}{P} - 1, & 0 < P < \gamma, \\
    0, & P \ge \gamma.
  \end{cases}
\end{equation*}

The scenario MPC problem for energy management is
\begin{subequations}
\label{eq:energy-management-scenario-MPC}
\begin{align}
\minimize\; &{\frac{1}{S}\sum_{s=1}^{S}\sum_{k=0}^{N-1} c_k(q_{s,k},\, u_k,\, d_{s,k})}\\
   \text{s.t.}\; &q_{s,k+1} = q_{s,k} - \frac{\Delta t}{b}\,P^B_k,
    \label{eq:em-dynamics}\\
    &q_{s,0} = q_t,
    \; s = 1,\ldots,S,
    \label{eq:em-initial-condition}\\
    &B_{\min} \le P^B_k \le B_{\max},
    \; k = 0,\ldots,N-1,
    \label{eq:em-BESS-bounds}\\
    &Q_{\min} \le q_{s,k} \le Q_{\max},
    \label{eq:em-SOC-bounds}\\
    &Q_N \le q_{s,N},
    \; s = 1,\ldots,S,
    \label{eq:em-terminal}\\
    &P\mathbf{u} = 0,\label{eq:em-non-anticipativity}
\end{align}
\end{subequations}
where $B_{\max} > 0$ is the maximum discharge power,
$B_{\min} < 0$ is the maximum charge power,
$Q_{\min}$ and $Q_{\max}$ are the SOC limits,
and $Q_N$ is the minimum terminal SOC reserve.
The non-anticipativity constraint~\eqref{eq:em-non-anticipativity} enforces a scenario-independent schedule,
following the structure of~\eqref{eq:scenario-MPC:non-anticipativity}.

In order to apply the proposed framework for the formulation~\eqref{eq:energy-management-scenario-MPC}, we need to reformulate it to fit the general scenario MPC.
The stage cost~\eqref{eq:em-stage-cost} depends explicitly on the disturbance $d_{s,k}$ through $D_{s,k}$,
whereas the general scenario MPC~\eqref{eq:scenario-MPC} requires the stage cost $c_k(\cdot, \cdot)$ to depend only on the state and control.
This mismatch arises because the SOC dynamics~\eqref{eq:em-SOC-dynamics} are disturbance-free, \ie
the disturbance does not enter the state trajectory and therefore cannot be captured by the state alone.
We augment the state with the disturbance $d_{s,k}$, thereby moving it from the stage cost into the dynamics.

Define the augmented state $\tilde{x}_{s,k} = [q_{s,k}, d_{s,k}]^\top$,
the augmented dynamics then take the linear form
\begin{equation*}
\label{eq:em-augmented-dynamics}
  \tilde{x}_{s,k+1}
  =
  \underbrace{\mm{1, 0; 0, 0}}_{\tilde{A}}\,\tilde{x}_{s,k}
  +
  \underbrace{\mm{0, -\tfrac{\Delta t}{b}; 0, 0}}_{\tilde{B}}\,u_k
  +
  \underbrace{\mm{0; d_{s,k+1}}}_{\tilde{d}_{s,k}}\text.
\end{equation*}
The exogenous term $\tilde{d}_{s,k} = [0,\, d_{s,k+1}]^\top$ is the disturbance that enters the augmented dynamics at step $k$,
corresponding to $d_{s,t+k}$ in the general formulation~\eqref{eq:scenario-MPC:dynamics}.
With the augmented state, the stage cost~\eqref{eq:em-stage-cost} depends on $\tilde{x}_{s,k}$ and $u_{s,k}$ only:
\begin{multline}
  \label{eq:em-stage-cost-augmented}
  c_k(\tilde{x}_{s,k},\, u_{s,k})
  = p_k \Delta t \!\left(P^c_k + d_{s,k} - P^B_k + \frac{1-\mu}{2\sqrt{\mu}}\,|P^B_k|\right)\\
  + p_p \max\!\left(P^c_k + d_{s,k} - P^B_k,\, 0\right)
  + \eta\, g(P^c_k)\text,
\end{multline}
where $d_{s,k} = [\tilde{x}_{s,k}]_2$ is read from the second component of the augmented state.
This matches the form $c_k(x_{s,k}, u_{s,k})$ required by~\eqref{eq:scenario-MPC}.
The augmented initial condition is
\begin{equation}
\label{eq:em-augmented-initial}
  \tilde{x}_{s,0} = [q_t, d_{s,0}]^\top\text,
\end{equation}
where $q_t$ is the current SOC (scenario-independent) and $d_{s,0}$ is the initial disturbance realization for scenario $s$.
The scenario-dependent component $d_{s,0}$ can be viewed as arising from the exogenous injection $\tilde{d}_{s,-1} = [0,\, d_{s,0}]^\top$ applied to the common base state $\tilde{x}_t = [q_t,\, 0]^\top$ before the first horizon step.
This is the scenario-dependent initial condition $x_{s,0}$.

In summary, the energy management scenario MPC~\eqref{eq:energy-management-scenario-MPC} with the augmented state is an instance of the general scenario MPC~\eqref{eq:scenario-MPC} with 
$A = \tilde{A},$
$B = \tilde{B}$,
$d_{s,t+k} = \tilde{d}_{s,k},$
$x_{s,0} = \tilde{x}_{s,0}$,
$x_{s,k} = \tilde{x}_{s,k}$,
cost function~\eqref{eq:em-stage-cost-augmented},
input constraint~\eqref{eq:em-BESS-bounds},
state constraints~\eqref{eq:em-SOC-bounds}--\eqref{eq:em-terminal},
and non-anticipativity constraint~\eqref{eq:em-non-anticipativity}.

\subsection{Experiment setting}
The problem parameters are set based on~\cite{cortes-aguirreEconomicMPC2025} as:
$\Delta t = \qty{0.25}{\hour}$, $p_k = \$0.1/\unit{\kWh}$, $p_p = \$19.19/\unit{\kW}$, $\eta = \$10$, $q_{s,0} = 0.5$, $B_{\min} = \qty{-700}{\kW}$, $B_{\max} = \qty{700}{\kW}$, $b = \qty{500}{\kWh}$, $\mu = 0.8$, $Q_{\min} = 0.2$, $Q_{\max} = 0.8$, and $\gamma = \qty{50}{\kW}$.
We perform a 10-step closed-loop MPC simulation with horizon $N \in \{48, 96\}$ and scenario $S \in \{10, 100\}$. 
Table \ref{tab:NN-hyperparameters} summarizes the training and architecture hyperparameters used for the neural network.

To evaluate the performance of the proposed learning-accelerated \SBMPC method, abbreviated as LA-\SBMPC, 
and compare it against IPOPT and MadNLP~\cite{shin2021graph}, we implement all algorithms in Julia on an Apple M4 Pro chip with 24 GB of RAM.

\begin{table}[!tb]
\centering
\caption{Training configuration and architecture parameters of the input-convex neural network.}
\label{tab:NN-hyperparameters}
\scalebox{1.0}{
\begin{tabular}{@{}ll ll@{}}
\toprule
Setting & Value & Setting & Value \\
\midrule
Optimizer              & AdamW     & Train/test instances         & 8000 / 2000             \\
Learning rate          & $10^{-3}$ & Hidden layers & $16 \times 16$ \\
Epochs                 & 5000      & Activation         & Softplus       \\
\bottomrule
\end{tabular}%
}
\end{table}

\subsection{Scenario construction and sampling}
Recall that the uncertainties in~\eqref{eq:energy-management-scenario-MPC} are the uncontrollable load $P^u$ and photovoltaic (PV) renewable generation $P^R$.
Following the scenario MPC formulation in~\eqref{eq:scenario-MPC}, each scenario must provide one realization of these uncertain signals over the prediction horizon.
We use the 24-hour real-world data from~\cite{cortes-aguirreEconomicMPC2025}, sampled every $\Delta t=0.25$ h, and repeat the profile over 7 days with additive Gaussian perturbations of standard deviation $1$ kW to construct the forecasting data.
This produces a scenario set with $L=7\times 96$ time steps.

For each signal type $\zeta\in\{u,R\}$, a Gaussian process (GP) model~\cite{williams2006gaussian} is fitted using time as the input, a constant mean function, and a periodic kernel with a 24-hour period.
The fitted GP provides the joint predictive distribution over the full
scenario horizon,
$$\tilde{\mathbf P}^{\zeta}_{s}
\sim
\mathcal{N}\left(\boldsymbol{\mu}^{\zeta},
\boldsymbol{\Sigma}^{\zeta}\right),\qquad s=1,\ldots,S,$$
where $\boldsymbol{\mu}^{\zeta}\in\mathbb{R}^{L}$ and
$\boldsymbol{\Sigma}^{\zeta}\in\mathbb{R}^{L\times L}$ are the GP
posterior mean and covariance evaluated at the 15-minute time grid.
The final scenario
values are projected onto the nonnegative range
$
\label{eq:gp-scenario-sampling}
P^\zeta_{s,k}=\max\{0,\tilde P^\zeta_{s,k}\}\text,
$
to enforce nonnegative load and generation values.

At each MPC sampling time $t$, we randomly select $S$ trajectories from this scenario set and extract the prediction window from $t$ to $t+N-1$.
Thus, the load and PV used by the MPC at a given start time $t$ are
\begin{equation*}
\mathbf{P}^{u}_{t} = P^u_{i,\,t:t+N-1}, \qquad \mathbf{P}^{R}_{t} = P^R_{i,\,t:t+N-1},
\end{equation*}
where $i$ denotes the randomly selected scenario indices.
Therefore, the length of each scenario window is determined by the MPC horizon $N$.
In the simulations, for each sampled scenario and each MPC window, the disturbance sequence is obtained from these two signals using~\eqref{eq:em-disturbance} and is passed to the optimization model through the power balance constraint~\eqref{eq:em-grid-exchange}.

\subsection{Results}
\begin{table}[!tb]
\centering
\caption{Average per-step closed-loop solving time (milliseconds) comparison.}
\label{tab:solving-time-comparison}
\setlength{\aboverulesep}{0pt}
\setlength{\belowrulesep}{0pt}
\renewcommand{\arraystretch}{1.3} 

\scalebox{1.0}{\begin{tabular}{lccc}
\toprule
Scenario & IPOPT & MadNLP & LA-SBMPC \\ \midrule
S = 10, N = 48 & 202 & 226 & \cellcolor{orange!30}{42} \\ \midrule
S = 100, N = 48 & 1879 & 6046 & \cellcolor{orange!30}{355} \\ \midrule
S = 10, N = 96 & 439 & 713 & \cellcolor{orange!30}{92} \\ \midrule
S = 100, N = 96 & 4632 & 11216 & \cellcolor{orange!30}{548} \\ \bottomrule
\end{tabular}%
}
\end{table}

\begin{figure}[!tb]
    \centering
    \includegraphics[width=0.9\linewidth]{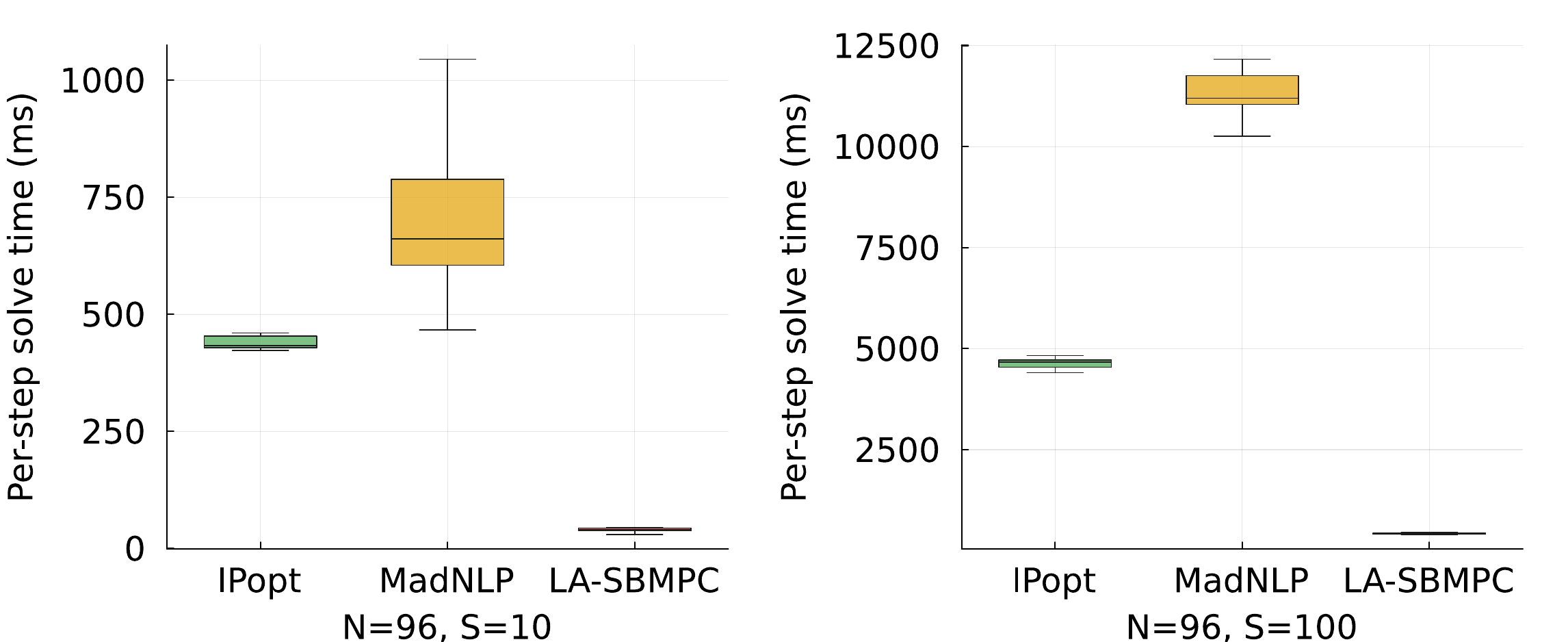}
    \caption{Per-step closed-loop solving time comparison.}     \label{fig:run-time-per-step}
\end{figure}

Table~\ref{tab:solving-time-comparison} compares the average per-step closed-loop solving time under different numbers of scenarios and prediction horizons. In the smallest setting with $S = 10$, $N = 48$, the proposed LA-\SBMPC  requires only $42$ ms, which is about $4.8\times$ faster than IPOPT and $5.4\times$ faster than MadNLP.
When the prediction horizon is increased to $N = 96$ while keeping $S = 10$, LA-\SBMPC  maintains its advantage, requiring $92$ ms and achieving speedups of $4.8\times$ and $7.8\times$ over IPOPT and MadNLP, respectively.
As the number of scenarios scales up to $S = 100$, the computational benefits of LA-\SBMPC  become even more significant.
Across the $S = 100$ settings, LA-\SBMPC achieves speedups ranging from $5.3\times$ to $8.5\times$ against IPOPT, and from $17\times$ to $20.5\times$ against MadNLP, with a maximum average per-step solve time of just $548$ ms.
These results suggest that LA-\SBMPC  scales well as the number of scenarios and the prediction horizon increase. 
Moreover, Fig.~\ref{fig:run-time-per-step} provides a per-step view of the closed-loop solving time, showing results consistent with Table~\ref{tab:solving-time-comparison}.
Specifically, LA-\SBMPC  achieves a clear runtime advantage in the larger settings.
As the number of scenarios and the prediction horizon increase, the baseline solvers (IPOPT and MadNLP) require substantially more time at each MPC step, whereas LA-\SBMPC  grows more moderately. 

\begin{figure}[!tb]
    \centering
    \includegraphics[width=0.9\linewidth]{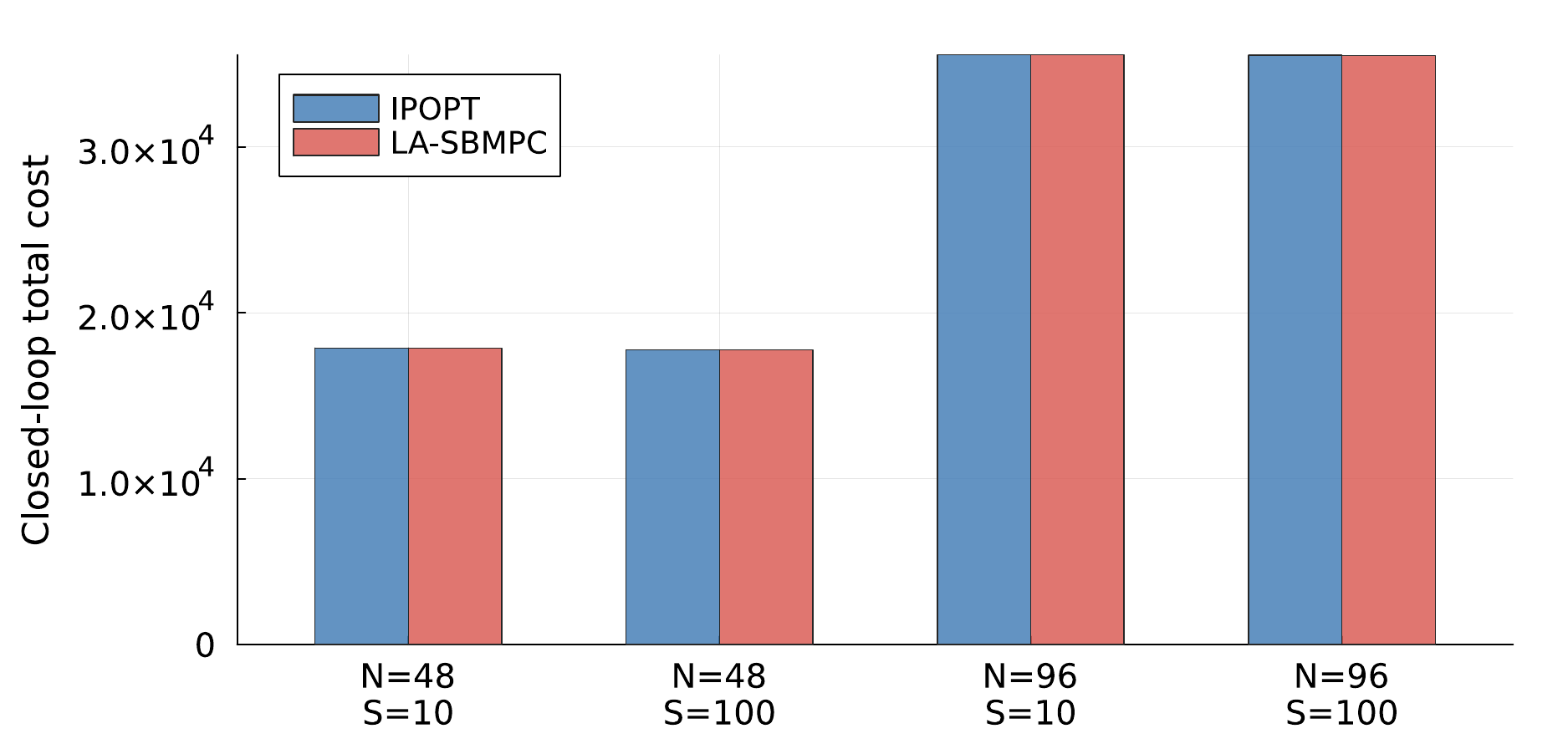}
    \caption{Closed-loop total cost comparison under varying scenario and horizon settings.}    
    \label{fig:total-cost}
\end{figure}

\begin{figure}[!tb]
    \centering
    \includegraphics[width=0.95\linewidth]{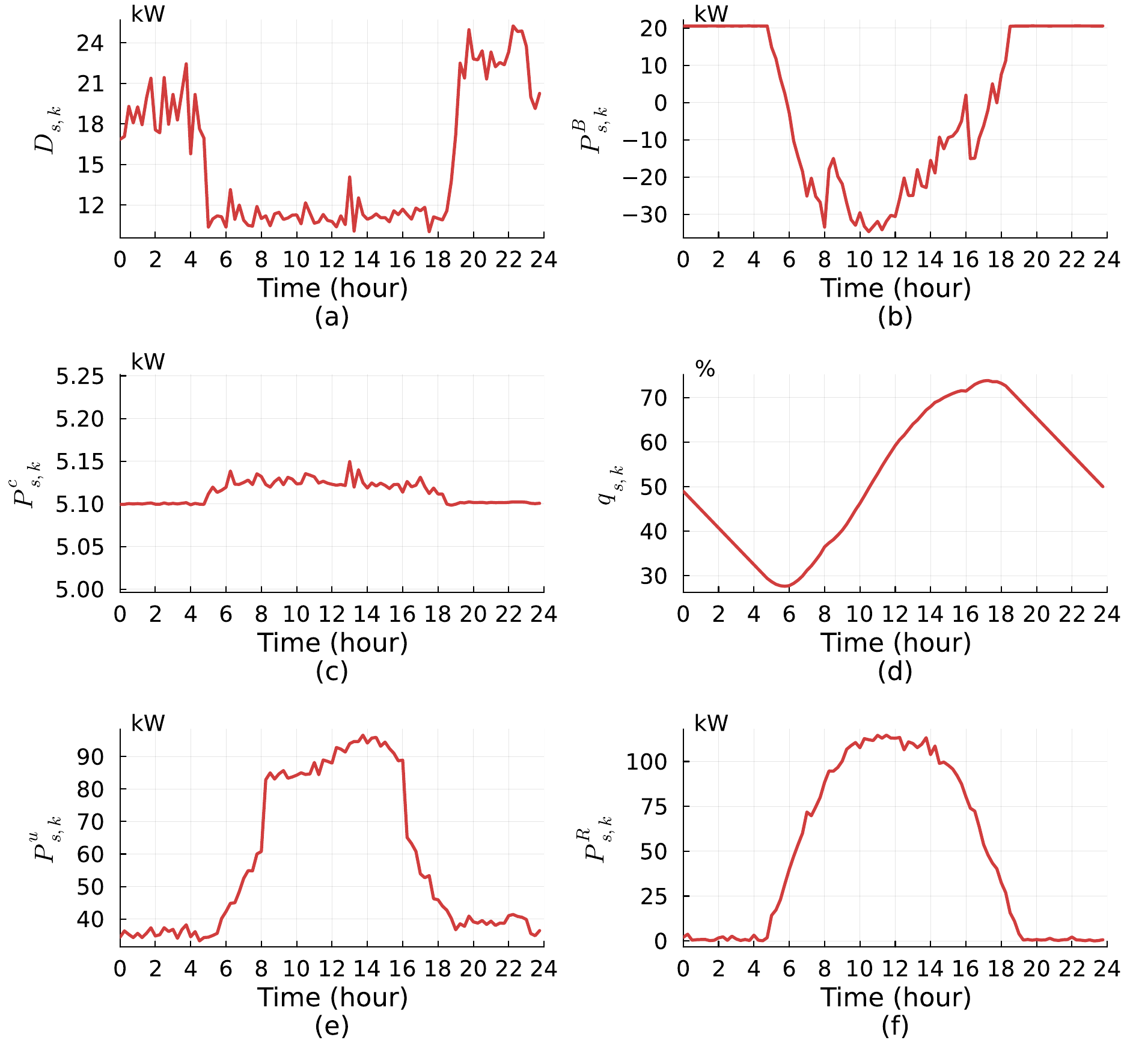}
    \caption{Closed-loop 24-hour operation of the microgrid scheduled by LA-\SBMPC  for $S=10$.}     \label{fig:LA-sMPC-prediction}
\end{figure}

Table~\ref{tab:constraint-violation} and
Fig.~\ref{fig:total-cost} evaluate the accuracy and feasibility of LA-SBMPC. The mean optimality gap is below $0.017\%$ in all cases, and the maximum gap is below $0.081\%$, showing close agreement with IPOPT.
The SOC, terminal, and non-anticipativity residuals are all small, confirming that LA-SBMPC maintains feasibility across the tested horizons and scenario counts.
Fig.~\ref{fig:LA-sMPC-prediction} illustrates the closed-loop performance of the microgrid over a 24-hour period under the LA-\SBMPC  schedule.
The figure includes: (a) the grid power exchange, (b) the BESS charging/discharging power, (c) the controllable load, (d) the BESS state of charge, (e) the uncontrollable load, and (f) the renewable generation. 
The resulting trajectories are physically reasonable and capture the main operating characteristics of the microgrid. 
In particular, the BESS charges or discharges in response to variations in load and renewable generation, while the SOC remains within its prescribed bounds. 
The controllable load is also adjusted according to system conditions, and the grid exchange reflects the balance between demand, renewable supply, and battery operation. 
These results demonstrate that the proposed method can produce meaningful closed-loop control actions for microgrid energy management over a period of time.

\begin{table}[!tb]
\centering
\caption{Optimality gaps and constraint residuals of LA-SBMPC over all benchmark scenarios.}
\label{tab:constraint-violation}
\setlength{\aboverulesep}{0pt}
\setlength{\belowrulesep}{0pt}
\renewcommand{\arraystretch}{1.15}
\scriptsize
\resizebox{\linewidth}{!}{\begin{tabular}{ccccc}
\toprule
$(N,S)$ & \begin{tabular}{c} Optimalitya gap (\%) \\ mean (max) \end{tabular} & \begin{tabular}{c} SOC \\ residual \end{tabular} & \begin{tabular}{c} Terminal \\ residual \end{tabular} & \begin{tabular}{c} Non-anticip. \\ residual \end{tabular} \\
\midrule
(48, 10)  & $2.40\times10^{-3}$ ($5.00\times10^{-3}$)  & $1.48\times10^{-13}$ & $1.67\times10^{-15}$ & $2.17\times10^{-15}$ \\
\midrule
(48, 100) & $6.00\times10^{-3}$ ($9.40\times10^{-3}$) & $4.83\times10^{-15}$ & $1.11\times10^{-16}$ & $1.32\times10^{-14}$ \\
\midrule
(96, 10)  & $1.56\times10^{-2}$ ($8.00\times10^{-2}$) & $8.28\times10^{-13}$  & $7.62\times10^{-14}$ & $7.11\times10^{-15}$ \\
\midrule
(96, 100) & $1.61\times10^{-2}$ ($8.04\times10^{-2}$) & $1.30\times10^{-12}$  & $8.30\times10^{-14}$ & $1.42\times10^{-14}$ \\
\bottomrule
\end{tabular}%
}
\end{table}

\section{Conclusion}\label{sec:conclusion}
In this paper, we accelerated the solution of computationally demanding scenario-based MPC problems by developing a learning-based ADMM framework.
By reformulating the scenario optimal control problem into a consensus form, we decoupled the coupling constraints and decomposed the problem via ADMM.
To overcome the latency of repeated primal updates within this structure, we integrated an input-convex neural network architecture to approximate the Moreau envelope of the local cost functions, thereby reducing the per-iteration computational overhead while strictly preserving the convexity of the problem.
The efficacy of this approach was demonstrated on a microgrid energy scheduling problem subject to load and renewable generation uncertainties. 
Our empirical results show that the proposed algorithm achieves runtime speedups compared to popular nonlinear programming solvers. 
Future work will focus on learning warm-start solutions and considering closed-loop stability guarantees.

\bibliographystyle{IEEEtran}
\bibliography{references}
\end{document}

%% file: intro.tex
\section{Introduction}\label{sec:introduction}
Model Predictive Control (MPC) is a popular framework for constrained optimal control, which computes feedback actions by repeatedly solving a finite-horizon optimal control problem~\cite{garcia1989model,nguyen2026ad}.
In many applications, however, the future evolution of the system is affected by uncertain disturbances.
These disturbances enter the predicted dynamics, constraints, and costs, and therefore must be accounted for when the controller plans actions over the horizon.
This issue is especially important in applications such as microgrid energy management, where load demand and renewable generation fluctuate over time and directly affect the economic and operational performance of the closed-loop system~\cite{zhao2023robust}.

Robust and stochastic MPC provide two common ways to handle uncertainty in finite-horizon control.
Robust MPC provides strong constraint-satisfaction guarantees by optimizing against worst-case uncertainty realizations, but often yields conservative and economically inefficient control actions~\cite{bemporad2007robust}.
Stochastic MPC instead models uncertainty probabilistically and typically enforces constraints in expectation or through chance constraints~\cite{mesbah2016stochastic}.
This probabilistic treatment can reduce conservatism, but the resulting chance-constrained problems can be nonconvex and computationally difficult to solve in real time~\cite{shapiro2021lectures}.

Scenario-based MPC (\SBMPC) offers a tractable deterministic approximation of stochastic MPC using sampled disturbance trajectories~\cite{micheli2022scenario}.
The scenario approach accommodates general disturbance distributions via sampling, but the problem grows quickly with both the prediction horizon and the number of scenarios.
To handle this growing complexity, the resulting scenario programs are often solved using the Alternating Direction Method of Multipliers (ADMM)~\cite{kang2015decomposition}.
This method is particularly attractive because it can exploit the separable scenario structure by duplicating variables and enforcing non-anticipativity through consensus constraints.
However, solving these primal updates across all scenarios at each iteration in real time remains a major computational bottleneck.

Learning-to-optimize has emerged as a promising tool for reducing the cost of repeated optimization in real-time decision-making~\cite{chen2022learning}.
For example, Baker~\cite{baker2019learning} learned warm-start points for AC optimal power flow to improve convergence speed.
Misra \emph{et al.}~\cite{misra2022learning} learned a small collection of active constraint sets for large-scale parametric optimization, thereby reducing online computation while preserving constraint satisfaction.
Our previous work developed the Learning-Enabled ADMM Framework (LEAF)~\cite{nguyen2026leaflearningenabledadmmframework}, which accelerates ADMM by learning the Moreau envelope associated with the primal update.
By leveraging the Moreau envelope to smoothly approximate objectives and recover proximal updates, LEAF efficiently accelerates the repeated local minimizations in ADMM-based optimization.

In this paper, we develop a learning-accelerated ADMM method for solving \SBMPC problems. 
The key idea is to first reformulate the \SBMPC problem in a consensus form that decomposes the ADMM primal step into parallel subproblems across the scenario-time grid.
We then apply LEAF to approximate these local primal updates, replacing costly repeated computations with efficient neural network evaluations. 
%
Our main contributions 
are as follows.
\begin{enumerate}[label=(\roman*)]
    \item We formulate a consensus ADMM decomposition for \SBMPC that decouples scenario-dependent dynamics and non-anticipativity constraints.
    The resulting primal updates are separable across scenarios and prediction time steps, which enables parallel computation over the full scenario-time grid.
    \item We integrate the Moreau-envelope learning framework LEAF into the decomposed ADMM algorithm to accelerate the local primal update subproblems.
    The learned update can be reused across scenarios with the same local cost structure, avoiding scenario-specific retraining while reducing the dominant per-iteration computation.
    \item We evaluate the proposed learning-accelerated \SBMPC method on a microgrid energy management problem with uncertain load demand and renewable generation.
    Numerical results show substantial speedups over IPOPT and MadNLP~\cite{shin2021graph} while maintaining high solution quality and meaningful closed-loop control performance.
\end{enumerate}
